\documentclass[runningheads]{llncs}
\usepackage{graphicx}
\usepackage{geometry}
\usepackage{amssymb}
\usepackage{amsmath}
\usepackage{hyperref}
\usepackage{xcolor}
\usepackage{enumerate}
\usepackage{listings}
\usepackage{complexity}
\usepackage{blkarray}
\usepackage{lmodern}

\begin{document}
\title{Inapproximability of shortest paths on perfect matching polytopes\thanks{The second author was supported by an ETH Postodctoral Fellowship.}}
%
%
\author{Jean Cardinal\inst{1}\orcidID{0000-0002-2312-0967} \and
Raphael Steiner\inst{2}\orcidID{0000-0002-4234-6136}}
\authorrunning{J. Cardinal and R. Steiner}
%
\institute{Universit\'e libre de Bruxelles (ULB), Brussels, Belgium\\
  \email{jean.cardinal@ulb.be}\and
  ETH Zurich, Z\"{u}rich, Switzerland\\
\email{raphaelmario.steiner@inf.ethz.ch}}
\maketitle              
\begin{abstract}
We consider the computational problem of finding short paths in the skeleton of the perfect matching polytope of a bipartite graph. We prove that unless $\P=\NP$, there is no polynomial-time algorithm that computes a path of constant length between two vertices at distance two of the perfect matching polytope of a bipartite graph. Conditioned on $\P\neq\NP$, this disproves a conjecture by Ito, Kakimura, Kamiyama, Kobayashi and Okamoto [SIAM Journal on Discrete Mathematics, 36(2), pp. 1102-1123 (2022)]. Assuming the Exponential Time Hypothesis we prove the stronger result that there exists no polynomial-time algorithm computing a path of length at most $\left(\frac{1}{4}-o(1)\right)\frac{\log N}{\log \log N}$ between two vertices at distance two of the perfect matching polytope of an $N$-vertex bipartite graph. These results remain true if the bipartite graph is restricted to be of maximum degree three. 

The above has the following interesting implication for the performance of pivot rules for the simplex algorithm on simply-structured combinatorial polytopes: If $\P\neq\NP$, then for every simplex pivot rule executable in polynomial time and every constant $k \in \mathbb{N}$ there exists a linear program on a perfect matching polytope and a starting vertex of the polytope such that the optimal solution can be reached in two monotone steps from the starting vertex, yet the pivot rule will require at least $k$ steps to reach the optimal solution. This result remains true in the more general setting of pivot rules for so-called \emph{circuit-augmentation algorithms}.

\end{abstract}

\section{Introduction}
The history of linear programming is intimately intertwined with that of Dantzig's simplex~algorithm.
While the simplex and its many variants are among the most studied algorithms ever,
a number of fundamental questions remain open. It is not known, for instance, whether there exists
a pivot rule that makes the simplex method run in strongly polynomial time.
Since the publication of the first examples of linear programs that make the original simplex algorithm run in exponential time, many alternative pivot rules have been proposed, fostering a tremendous amount of work in the past 75 years, both from the combinatorial and complexity-theoretic point of views.

The simplex algorithm follows a monotone path on the skeleton of the polytope defining the linear program. The following natural question was recently raised by De Loera, Kafer, and Sanit\`a~\cite{DKS22}:\\

\noindent\textit{``Can one hope to find a pivot rule that makes the simplex method use a shortest monotone path?''}.\\

As an answer, they proved that given an initial solution to a linear program, it is $\NP$-hard to find a 2-approximate shortest monotone path to an optimal solution.
It implies that unless $\P =\NP$, no polynomial-time pivot rule for the simplex can be guaranteed to reach an optimal solution in a minimum number of steps.

A similar result can also be deduced from two independent contributions, by Aichholzer, Cardinal, Huynh, Knauer, M\"utze, Steiner, and Vogtenhuber~\cite{ACHKMS21} on one hand, and by Ito, Kakimura, Kamiyama, Kobayashi, and Okamoto~\cite{IKKKO22} on the other hand. They proved that the above result holds for perfect matching polytopes of planar and bipartite graphs, albeit with a slightly weaker inapproximability factor of $3/2$ instead of $2$. Ito et al.~\cite{IKKKO22} conjecture that there exists a constant-factor approximation algorithm for the problem of finding a shortest path between two perfect matchings on the perfect matching polytope.

Our main result is a disproof of this conjecture under the $\P\not=\NP$ assumption: Strengthening the previous inapproximability results mentioned above, we show that unless $\P=\NP$ no $k$-approximation for a shortest path between two vertices at distance $2$ of a bipartite perfect matching polytope can be found in polynomial time, for any (arbitrarily large) choice of $k\in \mathbb{N}$. We also give an even stronger inapproximability result under the \emph{Exponential Time Hypothesis} (ETH).
The latter states that the 3-SAT problem cannot be solved in worst-case subexponential time, and is one of the main computational assumptions of the fine-grained complexity program~\cite{W18}.
As a consequence, there is not much hope of finding a pivot rule for the simplex algorithm yielding good approximations of the shortest path towards an optimal solution, even when the linear program is integer and its associated matrix totally unimodular.

\subsection{Our result}

We consider the complexity of computing short paths on the 0/1 polytope associated with perfect matchings of a bipartite graph.
Given a balanced bipartite graph $G=(V,E)$, where $V$ is partitioned into two equal-size independent sets $A$ and $B$, we define the
\emph{perfect matching polytope} $P_G\subseteq \mathbb{R}^E$ of $G$ as the convex hull of the 0/1 incidence vectors of perfect matchings of $G$.

It is well-known (see e.g. Chapter 18 in~\cite{schrijver}) that for \emph{bipartite} graphs $G$, there is a nice halfspace representation of $P_G$. An edge-vector $(x_e)_{e \in E} \in \mathbb{R}^E$ is in $P_G$ if and only if the following hold.

\begin{align}
    \sum_{e \ni v}{x_e} &= 1,  \qquad (\forall v \in V) \\ 
    x_e &\ge 0, \qquad (\forall e \in E).
\end{align}

The above is a compact encoding of $P_G$, with a number of constraints and variables of size polynomial in $G$. The assumption that $G$ is bipartite is crucial here: For non-bipartite $G$ the polytope defined by the above constraints has non-integral vertices and is thus not a representation of $P_G$~\cite{schrijver}. The matrix of this representation of a perfect matching polytope of a bipartite graph $G$ is simply the vertex-edge-incidence matrix of $G$, which is totally unimodular. The problem of maximizing a linear functional $w^T x$ subject to constraints $(1)$ and $(2)$ corresponds exactly to the problem of finding a perfect matching $M$ of $G$ whose weight $\sum_{e \in M}{w_e}$ is maximal. 

Given that the simplex algorithm moves along the edges of a polytope, it is crucial for our considerations to understand adjacency of vertices on $P_G$. The following result is well-known~\cite{C75,I02}.

\begin{lemma}\label{lemma:adjacency}
  For a bipartite graph $G$,
  two vertices of $P_G$ corresponding to two perfect matchings $M_1$ and $M_2$ are adjacent in the skeleton of $P_G$ if and only if the symmetric difference $M_1\Delta M_2$ is a cycle in $G$.
\end{lemma}

This cycle is said to be \emph{alternating} in both matchings, and one matching can be obtained from the other by \emph{flipping} this alternating cycle.
In general, we will say that two perfect matchings are \emph{at distance at most $k$} from each other on $P_G$, for some positive integer $k$, if one can be obtained from the other by successively flipping at most $k$ alternating cycles.

\begin{figure}
  \begin{center}
    \includegraphics[page=4, width=\textwidth]{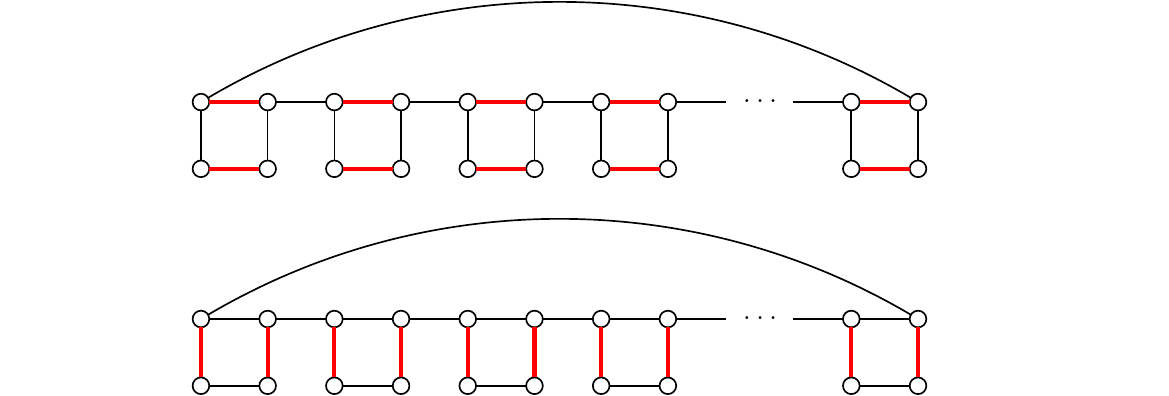}
  \end{center}
  \caption{\label{fig:ub}Two perfect matchings at distance two on the perfect matching polytope, but whose symmetric difference consists of an arbitrarily large number of even cycles.} 
  \end{figure}

Note that given any two perfect matchings $M_1$ and $M_2$ of a bipartite graph $G$, it is always the case that $M_1\Delta M_2$ is a collection of vertex-disjoint even cycles that are alternating in both matchings. The number of such cycles is therefore an upper bound on the distance between $M_1$ and $M_2$ on $P_G$. Interestingly, this upper bound can be arbitrarily larger than the actual distance. Figure~\ref{fig:ub} shows a construction of a graph $G$ with two matchings at distance two on $P_G$, whose symmetric difference consists of an arbitrary number of cycles.

Our main result is the following.

\begin{theorem}
  \label{thm:approxdistanceishard}
  Let $k \ge 2$ be any fixed integer.
  Unless $\P=\NP$, there does not exist any polynomial-time algorithm solving the following problem: 
  \begin{description}
  \item{\textbf{Input:}} A bipartite graph $G$ of maximum degree $3$ and a pair of perfect matchings $M_1, M_2$ of $G$ at distance at most $2$ on the polytope $P_G$.
  \item{\textbf{Output:}} A path from $M_1$ to $M_2$ in the skeleton of $P_G$, of length at most $k$.
  \end{description}

More strongly, for every absolute constant $\delta>0$, unless the Exponential Time Hypothesis fails, no polynomial-time algorithm can solve the above problem when $k$ is allowed to grow with the number $N$ of vertices of $G$ as $k(N)=\left\lfloor \left(\frac{1}{4}-\delta\right)\frac{\log N}{\log \log N} \right\rfloor$. 
\end{theorem}

A path on the perfect matching polytope of a bipartite graph $G$ is said to be \emph{monotone} with respect to some weight vector $w=(w_e)_{e \in E} \in \mathbb{R}^E$ on the edges of $G$ if the perfect matchings along the path have monotonically increasing total weights. Given two perfect matchings $M_1$ and $M_2$ at distance two on the perfect matching polytope, one can assign weights to edges so that (i) the path of length two between them is strictly monotone, and (ii) $M_2$ is the unique matching of maximal weight (this will be formally proven later in Lemma~\ref{lem:monotone}). This allows us to formulate our result as one about the hardness of reaching an optimal solution from a given feasible solution of a linear program on the perfect matching polytope.

\begin{corollary}
  \label{cor:monotone}
  Unless $\P=\NP$, there does not exist any polynomial-time constant-factor approximation algorithm for the following optimization problem: 
\begin{description}
\item{\textbf{Input:}} A bipartite graph $G=(V,E)$ of maximum degree $3$, a weight function $E\to\mathbb{R}^+$, and a perfect matching $M$ of $G$.
\item{\textbf{Output:}} A shortest monotone path on $P_G$ from $M$ to a maximum-weight perfect matching of $G$.
\end{description}

 Furthermore, assuming ETH, for an arbitrary but fixed $\delta>0$ no polynomial-time algorithm can achieve an approximation ratio of less than $\left(\frac{1}{8}-\delta\right)\frac{\log N}{\log \log N}$, where $N:=|V(G)|$. 
\end{corollary}

This corollary can be further interpreted as a statement on the existence of a polynomial-time pivot rule that would make the simplex method use an approximately shortest monotone path to a solution. Any such pivot rule could be used as an approximation algorithm for the above problem, contradicting the computational hypotheses.

\subsection{Pivot rules for circuit-augmentation algorithms.} 
Our work on distances in the skeleton of $P_G$ for bipartite graphs $G$ was originally motivated by questions regarding so-called \emph{circuit moves (or circuit augmentations)}, that have been recently studied in linear programming~\cite{BV22,DKS22,LHL15} as well as in the context of relaxations of the Hirsch conjecture concerning the diameter of polytopes~\cite{BFH15,KPS19}. A \emph{circuit move} extends the simplex-paradigm of moving along an incident edge of the constraint-polyhedron, by additionally allowing to move along certain non-edge directions, called \emph{circuits}. Given a linear program, the circuits in a well-defined sense represent all possible edge-directions that could occur after changing the right-hand side of the LP. The following is a formal definition. 

\begin{definition}[cf. Definition~1 in~\cite{DKS22}]\label{def:move}
Given a polyhedron of the form  $$\mathcal{P}=\{\mathbf{x}\in \mathbb{R}^n|A\mathbf{x}=\mathbf{b}, B\mathbf{x} \le \mathbf{d}\},$$ a \emph{circuit} is a vector $\mathbf{g} \in \mathbb{R}^n\setminus \{\mathbf{0}\}$ such that
\begin{enumerate}
    \item $A\mathbf{g}=0$, and
    \item the support of $B\mathbf{g}$ is inclusion-wise minimal among the collection $\{B\mathbf{y}|A\mathbf{y}=\mathbf{0}, y \neq \mathbf{0}\}$.
\end{enumerate}
\end{definition}

Given an LP $\{\max \mathbf{c}^T\mathbf{x}|\mathbf{x} \in \mathcal{P}\}$ for a polyhedron $\mathcal{P}$, a current feasible solution $\mathbf{x} \in \mathcal{P}$ and a circuit $\mathbf{g}$ with $\mathbf{c}^T\mathbf{g}>0$, a \emph{circuit move} then consists of moving to a new feasible solution $\mathbf{x}'=\mathbf{x}+t^\ast \mathbf{g}$, where $t^\ast \ge 0$ is maximal w.r.t. $\mathbf{x}+t^\ast \mathbf{g} \in \mathcal{P}$. Note that in general, an optimization algorithm based on a pivot rule for circuit moves may traverse several non-vertices of the polyhedron before reaching an optimal solution. 

Our interest in the perfect matching polytope for understanding the complexity of circuit-pivot algorithms came from the following observation. 

\begin{lemma}\label{lemma:circuits}
Let $G$ be a bipartite graph. Then 
\begin{enumerate}
    \item A non-zero vector $\mathbf{g} \in \mathbb{R}^E$ is a circuit for $P_G$ (defined by the (in)equalities $(1)$ and $(2)$) if and only if there exists an even cycle $C$ in $G$ decomposed into matchings $C^+, C^- \subseteq C$ (alternately appearing along $C$) and $\alpha>0$ such that $\mathbf{g}=\mathbf{g}(C,\alpha)$, where $$\mathbf{g}(C,\alpha)_e:=\begin{cases} \alpha, & \text{if }e \in C^+, \cr -\alpha, & \text{if }e \in C^-, \cr 0, & \text{if }e\notin C. \end{cases}.$$
    \item If $x$ is a vertex of $P_G$, and $x'\neq x$ is obtained from $x$ by a circuit move, then $x'$ is also a vertex of $P_G$ and adjacent to $x$ on the skeleton of $P_G$. 
\end{enumerate}
\end{lemma}

This observation implies that any sequence of circuit moves, applied starting from a vertex of $P_G$, will follow a monotone path on the skeleton of $P_G$ from vertex to vertex. Consequently, Corollary~\ref{cor:monotone} also yields an inapproximability result for polynomial pivot rules for circuit augmentation, as follows. 
\begin{corollary}\label{cor:circuitdistanceisalsohard}
  Unless $\P=\NP$, there does not exist a polynomial-time constant-factor approximation algorithm for the following problem. 
\begin{description}
\item{\textbf{Input:}} A bipartite graph $G$ of maximum degree $3$, a vertex $\mathbf{x} \in P_G$ and a linear objective function. 
\item{\textbf{Output:}} A shortest sequence of circuit moves on $P_G$ from $\mathbf{x}$ to an optimal solution. 
\end{description}
Furthermore, assuming ETH, no polynomial-time alorithm can achieve an approximation ratio of less than $\left(\frac{1}{8}-\delta\right)\frac{\log N}{\log \log N}$, where $N:=|V(G)|$ and $\delta>0$ is a constant.
\end{corollary}

\subsection{Related works}

Our work relates to two main threads of research in combinatorics and computer science: one obviously related to the complexity of the simplex method and linear programming in general, and another more recent one, aiming at building a thorough understanding of the computational complexity of so-called \emph{combinatorial reconfiguration problems}.

\subsubsection*{Complexity of the simplex method.}

In 1972, Klee and Minty showed that the original simplex method had an exponential worst-case behavior on what came to be known as {Klee-Minty cubes}~\cite{KM72}.
Since then, many other variants have been shown to have exponential or superpolynomial lower bounds~\cite{J73a,B77,GS79,AF17,DFH19}, although subexponential rules are known~\cite{HZ15}.
More dramatic complexity-theoretic results have been obtained recently~\cite{APR14,DS19}. In particular, it was shown by Fearney and Savani~\cite{FS15} that Dantzig's original simplex method can solve \PSPACE-complete problems: Given an initial vertex, deciding whether some variable will ever be chosen by the algorithm to enter the basis is \PSPACE-complete. The simplex method is also a key motivation for studying the diameter of polytopes, in particular the Hirsch conjecture, refuted in 2012 by Santos~\cite{S12}.

The hardness result on approximating monotone paths given by De Loera, Kafer, and Sanit\`a~\cite{DKS22} is in fact a corollary of the \NP-hardness of the following problem:
Given a feasible extreme point solution of the bipartite matching polytope and an objective function, decide whether there is a neighbor extreme point that is optimal.
A related result for circulation polytopes was proved by Barahona and Tardos~\cite{BT89}.
These two results, as well as the hardness results from Aichholzer et al.~\cite{ACHKMS21} and Ito et al.~\cite{IKKKO22} rely on the \NP-hardness of the Hamiltonian cycle problem.
In order to deal with the approximability of the shortest path, we have to resort to more recent inapproximability results on the longest cycle problem~\cite{BHK04}.

\subsubsection*{Reconfiguration of matchings.}

The field of \emph{combinatorial reconfiguration} deals with the problems of transforming a given discrete structure, typically a feasible solution of a combinatorial optimization problem, into another one using elementary combinatorial moves~\cite{IDHPSUU11,H13,N18,GIKO22}. The reachability problem, for instance, asks whether there exists such a transformation, while the shortest reconfiguration path problem asks for the minimum number of elementary moves.

A number of recent works in this vein deal with reconfiguration of matchings in graphs~\cite{IDHPSUU11,KMM12,BBHIKMMW19,BHIM19,GKM19}. Ito et al.~\cite{IDHPSUU11} proved that the reachability problem between matchings of size at least some input number $k$ and under single edge addition or removal was solvable in polynomial time. This was extended to an adjacency relation involving two edges by Kaminsk\'i et al.~\cite{KMM12}. The problem of finding the shortest reconfiguration path under this model was shown to be \NP-hard~\cite{BHIM19,GKM19}. 
Another line of work involves flip graphs on perfect matchings in which the adjacency relation corresponds to flips of alternating cycles of length exactly four~\cite{DH98,DH02,MFHHWU09,BBHIKMMW19,CRT21}. Note that for bipartite graphs, this flip graph is precisely the subgraph of the skeleton of the perfect matching polytope that consists of edges of length two.
Bonamy et al.~\cite{BBHIKMMW19} proved that the reachability problem in these flip graphs is \PSPACE-complete.
 

\section{Proof of Theorem~\ref{thm:approxdistanceishard}}

\subsection{Preliminaries}

First note that perfect matchings of a balanced bipartite graph $G=(A\cup B,E)$ can also be represented by orientations of $G$ in which every vertex in $A$ has outdegree one and every vertex in $B$ has indegree one. The edges of the matching are those oriented from $A$ to $B$. Alternating cycles in a perfect matching are one-to-one with directed cycles in this orientation, and flipping the cycle amounts to reverting the orientations of all its arcs. We will switch from one representation to another when convenient. 

We prove Theorem~\ref{thm:approxdistanceishard} by reducing from the problem of approximating the longest directed cycle in a digraph. We rely on the following two results from Bj\"{o}rklund, Husfeldt, and Khanna given as Theorems 1 and 2 in~\cite{BHK04}. 

\begin{theorem}[Bj\"{o}rklund, Husfeldt, Khanna~\cite{BHK04}]\label{thm:longcycle}
  Consider the problem of computing a long directed cycle in a given Hamiltonian digraph $D$ on $n$ vertices.
  \begin{enumerate}
  \item For every fixed $\varepsilon>0$, unless $\P=\NP$, there does not exist any polynomial-time algorithm that returns a directed cycle of length at least $n^{\varepsilon}$ in $D$.
  \item For every polynomial-time computable increasing function $f:\mathbb{N}\to \mathbb{N}$ in $\omega(1)$, unless the Exponential Time Hypothesis fails, there does not exist any
    polynomial-time algorithm that returns a directed cycle of length at least $f(n) \log n$ in $D$.
    \end{enumerate}
\end{theorem}

Note that in the two problems, the input graph is guaranteed to be Hamiltonian, yet it remains hard to explicitly \emph{construct} a directed cycle of some guaranteed length. Characterising the approximability of the longest cycle problem in undirected graphs is a longstanding open question~\cite{AYZ95,GN08}. 

The second ingredient of our proof is the following lemma, perhaps of independent interest, that bounds the increase in length of a longest directed cycle after a number of cycle flips in a digraph.

\begin{lemma}
  \label{lem:cycle}
  Let $G$ be an undirected graph, and let $C_1,\ldots,C_t$ be a sequence of (not necessarily distinct) cycles in $G$.
  Let $D_0, D_1, \ldots, D_t$ be a sequence of orientations of $G$ such that for each $i \in [t]$ the cycle $C_i$ is directed in $D_{i-1}$ and such that $D_i$ is obtained from $D_{i-1}$ by flipping $C_i$. 

There exists a polynomial-time algorithm that, given as input a number $\ell$, the orientations $D_0,\ldots,D_t$ and a directed cycle $C$ in $D_t$ of length $|C|>\ell^{t+1}$, computes a directed cycle in $D_0$ of length at least $\ell$. 
\end{lemma}

\begin{figure}
  \begin{center}
    \includegraphics[page=5, width=\textwidth]{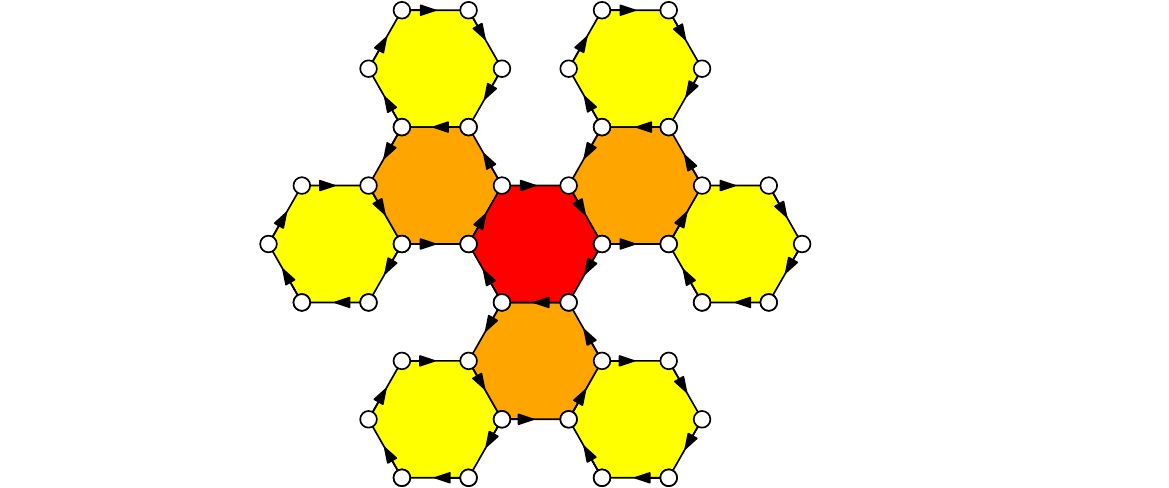}
  \end{center}
  \caption{\label{fig:cycle}A subcubic digraph in which the maximum length of a directed cycle is 6, yet after flipping two cycles, we reach an orientation of the same graph that has the outer face as a directed cycle, and thus contains a directed cycle of length greater than $3 \cdot 2 \cdot 5$ (exactly $42$ in this example). In general, by iterating the same construction, for every integer $k$ one can construct a digraph with maximum directed cycle length $\ell=2k$ and such that after flipping at most $t$ cycles, one can reach an orientation containing a directed cycle of length at least $k \cdot (k-1)^{t-1} \cdot (2k-1) \simeq 2\cdot (\ell / 2)^{t+1}$.} 
  \end{figure}

The bound of Lemma~\ref{lem:cycle} can be shown to be essentially tight. We refer to Figure~\ref{fig:cycle} for an illustration of the construction.
The proof of Lemma~\ref{lem:cycle} is deferred to the end of this section.

\subsection{Reduction}

\begin{figure}
  \begin{center}
    \includegraphics[page=1, width=\textwidth]{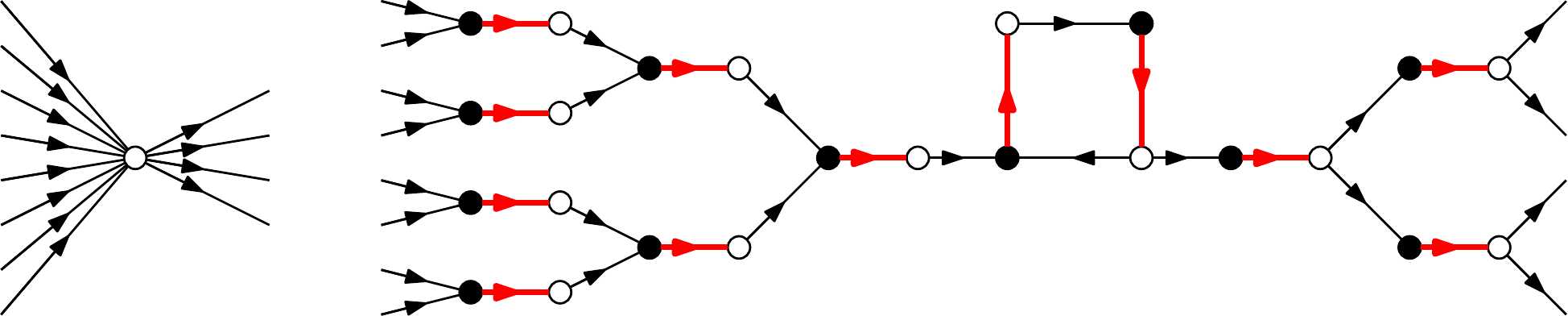}
  \end{center}
  \caption{\label{fig:reduction}Illustration of the reduction in the proof of Theorem~\ref{thm:approxdistanceishard}. Every vertex in the given Hamiltonian digraph $D$ (left) is replaced by the depicted gadget (right), yielding a maximum degree-three bipartite graph with a perfect matching.} 
  \end{figure}

We now give a proof of Theorem~\ref{thm:approxdistanceishard}, assuming Lemma~\ref{lem:cycle}.

\begin{proof}[Theorem~\ref{thm:approxdistanceishard}]
  We consider the first problem in Theorem~\ref{thm:longcycle}: For a fixed $\varepsilon>0$, given a Hamiltonian digraph $D$ on $n$ vertices, return a directed cycle of length at least $n^{\varepsilon}$.
  We first construct a digraph $D'$ from $D$ by replacing every vertex $v$ of $D$ by the gadget illustrated on Figure~\ref{fig:reduction}.
  The gadgets are obtained by applying the following transformations to every vertex $v$ of $D$:
  \begin{enumerate}
  \item The set of incoming arcs of $v$ is decomposed into a balanced binary tree with $\text{deg}_D^-(v)$ leaves and a degree-one root identified to $v$. Each internal node of this binary tree is further split into an arc. All arcs of the tree are oriented towards the root.
  \item The set of outgoing arcs are split into a tree with $\text{deg}_D^+(v)$ leaves in a similar fashion, with all arcs oriented away from the root.
  \item Finally, the vertex $v$ itself is replaced by a directed 4-cycle, such that the single incoming arc from the first tree and the single outgoing arc from the second tree have adjacent endpoints on the cycle.
  \end{enumerate}

  \begin{figure}
  \begin{center}
    \includegraphics[page=2, width=\textwidth]{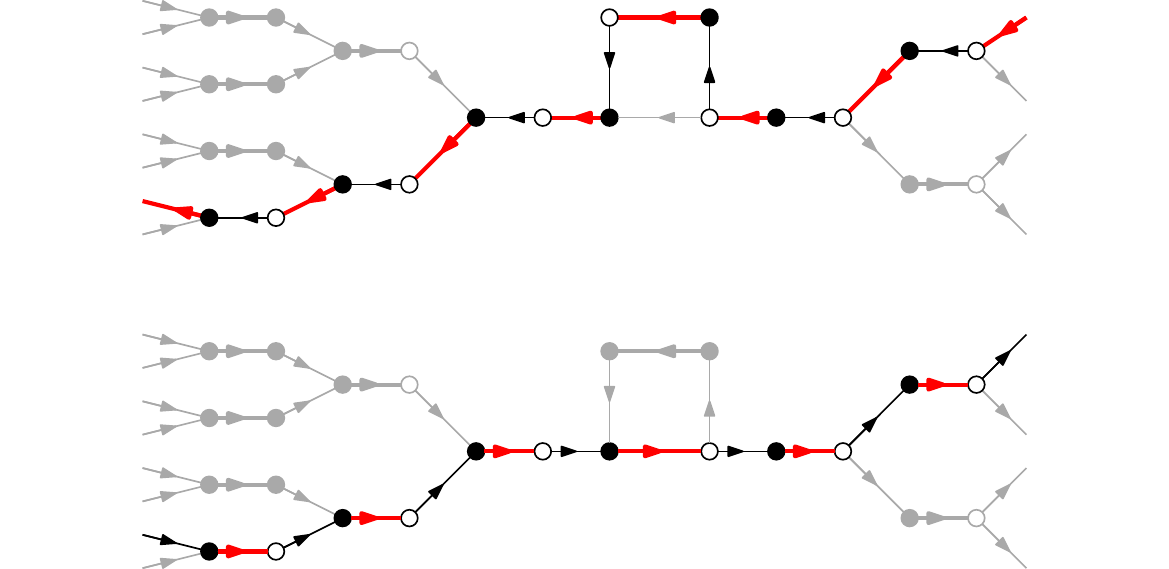}
  \end{center}
  \caption{\label{fig:flips}Flipping the 4-cycles of each gadget in $D'$ can be done with two successive cycle flips, using the Hamiltonian cycle of $D$.} 
  \end{figure}

  The digraph $D'$ thus obtained is bipartite and subcubic. Furthermore, it is easy to see by construction that for every vertex $v \in V(D)$, the corresponding gadget in $D'$ has at most $$4\text{deg}_D^-(v)+4+4\text{deg}_D^+(v) \le 8(n-1)+4<8n$$ vertices, such that $N:=|V(D')|< n\cdot 8n = 8n^2$, and $D'$ is of polynomial size.
  
  Furthermore, the orientation of $D'$ is such that every vertex in one side of the bipartition has outdegree one, and every vertex in the other has indegree one, hence it corresponds to a perfect matching $M_1$.
  By flipping the alternating 4-cycle in each gadget, we obtain another perfect matching $M_2$.
  We observe that $M_2$ can be obtained from $M_1$ in two cycle flips, by using the Hamiltonian cycle of $D$ twice (see Figure~\ref{fig:flips}). Hence, while $M_1\Delta M_2$ consists of $n$ disjoint 4-cycles, $M_2$ is in fact at distance two from $M_1$ on the perfect matching polytope of $D'$.
  The underlying undirected graph of $D'$ together with the two perfect matchings $M_1$ and $M_2$ therefore constitute an instance of the problem described in Theorem~\ref{thm:approxdistanceishard}.
  We now show that any sequence of length at most $k$ of alternating cycle flips transforming $M_1$ into $M_2$ can be turned in polynomial time into a cycle of length at least $n^{\varepsilon}$ in $D$, for some $\varepsilon>0$ depending solely on $k$.

  Consider a sequence of $k$ cycles $C'_1,C'_2,\ldots ,C'_k$ such that $C'_1$ is alternating in $D'$, and $C'_i$ is alternating in the graph obtained from $D'$ after flipping the cycles $C'_1,C'_2,\ldots ,C'_{i-1}$, in this order, and such that flipping all $k$ cycles in sequence transforms $M_1$ into $M_2$.
  Observe that the sum of the lengths of the cycles in this sequence must be at least $n$, since all the orientations of the $4$-cycles in the $n$ different gadgets in $D'$ have to be flipped, and since every single cycle $C_i'$ can intersect at most $|C_i'|$ different gadget-$4$-cycles.
  Let $\ell = \lceil n^{1/(k+2)}\rceil$.
  We have
  $$
  \sum_{i=1}^k |C'_i|\geq n =(1-o(1))\ell^{k+2} > \sum_{i=1}^k \ell^{i+1},
  $$
  hence from the pigeonhole principle, at least one cycle $C'_i$ in the sequence has length $|C'_i|>\ell^{i+1}$.
  From Lemma~\ref{lem:cycle}, we can now compute in polynomial time a directed cycle in $D'$ of length at least $n^{\varepsilon'}$ for $\varepsilon' = 1/(k+2)$. Let us call this cycle $C'$.
  Note that for every gadget in $D'$ corresponding to a vertex $v$ of $D$, either $C'$ is vertex-disjoint from this gadget, or it traverses it via exactly one directed path, consisting of a leaf-to-root path in the in-tree, a directed path of length $7$ touching the $4$-cycle of the gadget, and then a root-to-leaf path in the out-tree. 
  From this it follows that by contracting the edges of the gadgets, the cycle $C'$ in $D'$ can be mapped to a cycle $C$ in $D$. Note that the in- and out-degree of a vertex in $D$ is at most $n-1$, thus all the in- and out-trees in $D'$ corresponding to the gadgets have depth at most $2\lceil\log_2 n\rceil$. Consequently, the length of $C'$ can be shrinked by at most a factor of $4\lceil\log_2 n\rceil+7$ by contracting the gadgets. In other words, we obtain a directed cycle $C$ in $D$ of length at least $n^{\varepsilon'}/(4\lceil\log_2 n\rceil+7) = n^{\varepsilon}$, for $\varepsilon = \varepsilon' - o(1)$.
  Hence if we can obtain in polynomial time a sequence of at most $k=O(1)$ flips transforming $M_1$ into $M_2$, we can also find a cycle of length at least $n^{\varepsilon}$ in $D$ for some fixed $\varepsilon>0$. This establishes the first statement of Theorem~\ref{thm:approxdistanceishard}. 

  It remains to prove the second statement. We consider the second problem in Theorem~\ref{thm:longcycle}, in which we seek a path of length at least $f(n)\log n$, for some computable function $f(n)=\omega (1)$. Suppose that for some $\delta>0$ there is a polynomial-time algorithm that can find a sequence of at most $k=k(N)=\left\lfloor\left(\frac{1}{4}-\delta\right) \frac{\log N}{\log\log N}\right\rfloor $ flips transforming $M_1$ into $M_2$. Note that  $$k+2 \le  \left(\frac{1}{4}-\delta\right) \frac{\log 8n^2}{\log\log 8n^2}+2 < \frac{1}{2}(1-\delta)\frac{\log n}{\log \log n}$$ for $n$ large enough (larger than some constant depending only on $\delta$). 
  
  Now, from the same reasoning as above, we can turn such an algorithm into a polynomial-time algorithm that finds a directed cycle in $D$ of length at least
  $$
  \frac{n^{1/(k+2)}}{4\lceil\log_2 n\rceil+7} > \frac{n^{2 \cdot \log \log n/((1-\delta)\log n)}}{O(\log n)} = \Omega\left(\frac{\log^{2/(1-\delta)}n}{\log n}\right)= f(n)\log n,
  $$
  for a computable function $f(n) = \Omega(\log^{2\delta/(1-\delta)} n)=\omega(1)$.
  This, from Theorem~\ref{thm:longcycle}, is impossible unless the Exponential Time Hypothesis fails.
  \qed\end{proof}

  In order to deduce Corollary~\ref{cor:monotone} from Theorem~\ref{thm:approxdistanceishard}, we need the following lemma.
  
  \begin{lemma}
    \label{lem:monotone}
    Given a bipartite graph $G=(V,E)$, let $M_1$ and $M_2$ be two perfect matchings in $G$ at distance two on the perfect matching polytope, hence such that $M_2=(M_1\Delta C_1)\Delta C_2$ for some pair $C_1,C_2$ of cycles in $G$, and such that $M_1 \Delta C_1=:M'$ is also a perfect matching. Then there exists a weight function $w:E\to\mathbb{R}^+$ such that
    \begin{enumerate}
    \item $M_2$ is the unique maximum-weight perfect matching of $G$,
    \item $w(M_1)<w(M')<w(M_2)$ (where $w(M)=\sum_{e\in M}w(e)$).
    \end{enumerate}
  In other words, there exists a linear program over the perfect matching polytope of $G$ such that the path $M_1, M_1\Delta C_1=M'=M_2 \Delta C_2,M_2$ is a strictly monotone path and $M_2$ is the unique optimum.
  \end{lemma}
  \begin{proof}
    We first claim that $M_1 \setminus (M_2 \cup M')\neq \emptyset$. Indeed, the set of edges $M_2 \cup M'$ forms a subgraph of $G$ whose components are the edges in $M_2 \setminus C_2=M' \setminus C_2$ and the cycle $C_2$. From this it is easy to see that if we were to have that $M_1 \subseteq M_2 \cup M'$, then necessarily $M_1\in \{M_2, M'\}$, a contradiction to the fact that $M_1$ and $M_2$ have distance exactly two. Now pick an edge $e_0 \in M_1 \setminus (M_2 \cup M')$. Define a weighting $w:E\rightarrow \mathbb{R}$ of the edges in $G$ as follows: $w(e_0):=0$, $w(e):=1$ for every $e \in E\setminus (\{e_0\}\cup M_2)$ and $w(e):=1+\frac{1}{n}$ for every $e \in M_2$, where $n$ denotes the number of vertices of $G$. It is now easily verified that $M_2$ is the unique maximum-weight perfect matching in $G$, and that $w(M_1) \le 0+(\frac{n}{2}-1)(1+\frac{1}{n})<\frac{n}{2}\le w(M') < w(M_2)$. This proves the claim.  \qed
  \end{proof}
  
\subsection{Proof of Lemma~\ref{lem:cycle}}

\begin{proof}[Lemma~\ref{lem:cycle}]
Let the orientations $D_0,D_1,\ldots,D_t$ of $G$ be given as input, together with a directed cycle $C$ in $D_t$ and a number $\ell \in \mathbb{N}$ such that $|C|>\ell^{t+1}$.

Our algorithm starts by computing the sequence of cycles $C_1,\ldots,C_t$ by determining for each $i \in [t]$ the set of edges with different orientation in $D_{i-1}$ and $D_i$. Next we compute in polynomial time the subgraph $H$ of $G$ which is the union of the cycles $C_1,\ldots,C_t$ in $G$. We in particular compute a list of the vertex sets of its connected components, which we call $Z_1,\ldots,Z_c$ for some number $c \ge 1$.

For later use, let us prove the following fact:

\medskip

\textbf{Claim ($\ast$).} For each $r \in [c]$ the induced subdigraph $D_0[Z_r]$ of $D_0$ is strongly connected. 
\begin{proof}[Claim ($\ast$)]
For each $i \in [t]$, let us denote by $H_i$ the subgraph of $G$ consisting of the vertices and edges on the cycles $C_{t-i+1}, \ldots, C_t$.  Then we have $H_t=H$. We will establish the claim by proving inductively that for every $i=1,2,\ldots,t$, every connected component of $H_i$ induces a strongly connected subdigraph of $D_{t-i}$.

In the case $i=1$, this holds easily: The only connected component of $H_1$ is the cycle $C_t$, and since it forms a directed cycle both in $D_{t-1}$ and $D_t$, its vertices induce a strongly connected subdigraph of $D_{t-1}$, as desired. 

For the inductive step suppose we have established for some $1 \le i<t$ that all components of $H_i$ induces strongly connected subgraphs of $D_{t-i}$.

Note that when we move from $H_i$ to $H_{i+1}$, the only change that happens is the addition of the vertices and edges on the cycle $C_{t-i}$. If $C_{t-i}$ is vertex-disjoint from $H_i$, then it simply forms a new component in $H_{i+1}$ and since it is a directed cycle in $D_{t-(i+1)}$, it induces a strongly connected subdigraph thereof, while also all other components still induce strongly connected subdigraphs, as their orientations are not affected when flipping $C_{t-i}$. 

Otherwise, $C_{t-i}$ together with all components of $H_i$ intersected by it creates a new component $Z$ of $H_{i+1}$. Since each component intersected by $C_{t-i}$ induces a strongly connected subdigraph of $D_{t-i}$ by our assumption, and since $C_{t-i}$ is directed in $D_{t-i}$, it is easy to see that also the new component $X$ induces a strongly connected subdigraph of $D_{t-i}$. The only change when moving from $D_{t-i}$ to $D_{t-(i+1)}$ is that the orientations of the arcs of $C_{t-i}$ are reverted. However, it is easily seen that this operation preserves the fact that $X$ induces a strongly connected subdigraph, since flipping $C$ maintains the pairwise reachability of vertices on $C$ in the subdigraph induced by $X$. 

All components of $H_i$ not intersected by $C_{t-i}$ remain components also in $H_{i+1}$ and have the same orientations in $D_{t-(i+1)}$ as in $D_{t-i}$. Hence, also they induce strongly connected subdigraphs of $D_{t-(i+1)}$, as desired. This concludes the proof by the principle of induction. 
\qed\end{proof}

Let $(x_0,x_1,\ldots,x_{k-1},x_k=x_0)$ be the cyclic list of vertices on the directed cycle $C$ in $D_t$, with edges oriented from $x_i$ to $x_{i+1}$ for all $i \in [k-1]$. By assumption on the input, we have $k=|C|>\ell^{t+1}$. 

We first check if $C$ is vertex-disjoint from $H$, in which case we may return $C$, which is then also a directed cycle in $D_0$ of length $k>\ell^{t+1}\geq\ell$, as desired. 

Otherwise, $C$ intersects some of the components of $H$. We then for each vertex $x_i \in V(C)$ compute a label $\text{lab}(x_i) \in [c+1]$, defined as $\text{lab}(x_i):=r$ if $x_i \in Z_r$ lies in the $r$-th component of $H$, and $\text{lab}(x_i):=c+1$ if $x_i$ is not a vertex of $H$. We next compute an auxiliary weighted directed multigraph $M$ on the vertex set $[c]$ as follows: For every maximal subsequence of $C$, of length at least two, of the form $x_i,x_{i+1},\ldots, x_j$ (addition to be understood modulo $k$) such that $\text{lab}(x_s)=c+1$ for all $s=i+1,\ldots,j-1$ (if any), we add an additional arc from $\text{lab}(x_i)$ to $\text{lab}(x_j)$ and give it weight $j-i$, the corresponding number of arcs in $C$. Note that the total arc weight in $M$ is exactly $|C|$, while the total number of arcs is exactly $|V(C) \cap V(H)| \le |V(H)|$. The construction of $M$ is illustrated in Figure~\ref{fig:M}.

\begin{figure}
  \begin{center}
    \includegraphics[page=3, width=\textwidth]{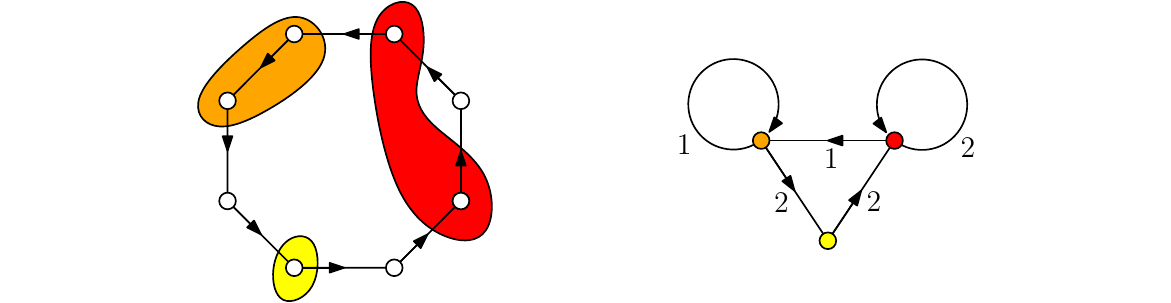}
  \end{center}
  \caption{\label{fig:M}Construction of the auxiliary directed multigraph $M$ in the proof of Lemma~\ref{lem:cycle}. The cycle $C$ is shown on the left, together with the connected components of $H$ that it intersects. The resulting weighted directed multigraph $M$ is shown on the right.}
  \end{figure}

Furthermore, by definition every vertex in $M$ has the same number of incoming and outgoing arcs. Hence, we may compute in polynomial time an edge-disjoint decomposition of $M$ into directed cycles (including possible loops) in $M$. Let $W_1,\ldots,W_p$ for some $p \in \mathbb{N}$ be the list of edge-disjoint directed cycles in this decomposition of $M$. We now create, for each $W_i$, a directed cycle $K_i$ in $D_0$ of length $|K_i|\ge \text{weight}(W_i)$, where $\text{weight}(W_i)$ is the total arc weight of $W_i$, as follows:

Let $(l_0,l_1,\ldots,l_s=l_0)$ be the cyclic vertex-sequence of $W_i$. 

For each arc $(l_j,l_{j+1})$ in $W_i$, we consider the corresponding subsequence $P(l_j,l_{j+1})$ of $C$ which starts in $Z_{l_j}$, ends in $Z_{l_{j+1}}$, and all whose internal vertices are not contained in $H$. We note that since arcs outside $H$ have the same orientation in $D_0$ and $D_t$, the subsequence $P(l_j,l_{j+1})$ is a directed path or a directed cycle also in $D_0$ which starts in $Z_{l_j}$ and ends in $Z_{l_{j+1}}$. 

We first check whether there exists an index $j$ such that $P(l_j,l_{j+1})$ is a directed cycle. In this case, necessarily $W_i$ is a loop (i.e. $s=0$) and $l_j=l_{j+1}=l_0$. We thus may simply put  $K_i:=P(l_j,l_{j+1})$, with $\text{weight}(W_i)=|K_i|$ satisfied by definition of the weights in $M$.

Otherwise, each of the $P(l_j,l_{j+1})$ is a directed path in $D_0$. We now make use of Claim~$(\ast)$, which tells us that $D_0[Z_{l_j}]$ is strongly connected for every $j=0,1,\ldots,s-1$. We may therefore compute in polynomial time for each $j=0,1,\ldots,s-1$ a directed path $Q_j$ in $D_0[Z_{l_j}]$ (possibly consisting of a single vertex) which connects the endpoint of $P(l_{j-1},l_j)$ to the starting point of $P(l_j,l_{j+1})$ (index-addition modulo $s$). Crucially, note that any two directed paths in the collection $\{P(l_j,l_{j+1}), Q_j| j=0,1,\ldots,s-1\}$ are vertex-disjoint except for shared common endpoints. We now compute the directed cycle $K_i$ in $D_0$, which is the union of the directed paths $P(l_j,l_{j+1})$ and the directed paths $Q_j$ for $j=0,\ldots,s-1$. It is clear that its length $|K_i|$ is lower-bounded by the sum of the lengths of the $P(l_j,l_{j+1})$, which by definition of $M$ equals the sum of arc-weights on $W_i$, i.e., we indeed have $|K_i| \ge \text{weight}(W_i)$ also in this case. 

After having computed the directed cycles $K_1,\ldots,K_p$ in $D_0$, the algorithm checks whether one of the cycles has length $|K_i|\ge \ell$.
If so, it returns the cycle $K_i$ and the algorithm stops with the desired output. Otherwise, we have $|K_i|<\ell$ for $i=1,\ldots,p$, which implies that
$$\ell^{t+1}<|C|=\text{weight}(M)=\sum_{i=1}^{p}{\text{weight}(W_i)} \le \sum_{i=1}^{p}{|K_i|} \le p(\ell -1).$$
Note that $p$ is at most as large as the number of arcs in $M$, which in turn is bounded by $|V(H)|$. We thus obtain
$$\ell^{t+1} < |V(H)| \cdot (\ell -1) \le \sum_{i=1}^{t} |C_i|\cdot (\ell -1).$$
This yields that
$$\sum_{i=1}^{t}{|C_i|} > \frac{\ell^{t+1}}{\ell -1}>\sum_{i=1}^{t}{\ell^i}.$$
Therefore there exists $i \in \{1,\ldots,t\}$ such that $|C_i|>\ell^i$. The algorithm proceeds by finding one cycle $C_i$ with this property. Note that $C_i$ is a directed cycle in the orientation $D_{i-1}$ of $G$. Hence a recursive call of the algorithm to the input $D_0,D_1,\ldots,D_{i-1}$ and the cycle $C_i$ will yield a directed cycle of length at least $\ell$ in $D_0$, as desired. 

This proves the correctness of the described algorithm. As all the steps between two recursive calls are executable in polynomial time in the size of $G$ and $t$, and since there will clearly be at most $t-1$ recursive calls in any execution of the algorithm, the whole algorithm runs in polynomial time, as desired. 
\qed\end{proof}

\subsection{Proof of Lemma~\ref{lemma:circuits}}
\begin{proof}[Lemma~\ref{lemma:circuits}]
\begin{enumerate}
    \item We start by noting that the hyperplane restrictions $(1)$ and $(2)$ for the perfect matching polytope $P_G$ can be expressed as $P_G=\{\mathbf{x} \in \mathbb{R}^E|A\mathbf{x}=\mathbf{b}, B\mathbf{x} \le \mathbf{d}\}$, where $A \in \mathbb{R}^{V \times E}$ is the vertex-edge incidence matrix of $G$, $\mathbf{b}=\mathbf{1}$, $B=-I_E$ and $\mathbf{d}=\mathbf{0}$. Applying Definition~\ref{def:move}, this means that a vector $\mathbf{g} \in \mathbb{R}^E\setminus\{\mathbf{0}\}$ is a circuit for $P_G$ if and only if it is a support-minimal member of the set $$\mathcal{Y}=\left\{\mathbf{y}\in \mathbb{R}^E\setminus\{\mathbf{0}\}\bigg\vert\sum_{e \ni v}{y_e}=0 \text{ for every } v \in V\right\}.$$
    Now consider a circuit $\mathbf{g} \in \mathcal{Y}$ and let $H$ be the subgraph of $G$ spanned by those $e \in E$ for which $\mathbf{g}_e \neq 0$. The defining equations of $\mathcal{Y}$ immediately imply that $H$ has minimum degree at least $2$. Thus $H$ contains a cycle $C$, which is necessarily even. It is easy to check that the vectors $\mathbf{g}(C,\alpha), \alpha > 0$ described in the lemma form members of $\mathcal{Y}$. The support-minimality of $\mathbf{g}$ thus implies that $\mathbf{g}$ is non-zero only on $C$. It now easily follows that $\mathbf{g}=\mathbf{g}(C,\alpha)$ for some $\alpha>0$. This concludes the proof of the first item. 
    \item Let $\mathbf{x}$ be a vertex of $P_G$ and $M$ the corresponding perfect matching of $G$. Let $\mathbf{x}' \in P_G$ be a vector reachable via a circuit-move from $\mathbf{x}$. This means that there is an even cycle $C$ in $G$ decomposed into two matchings $C^+$ and $C^-$, some $\alpha>0$ and $t^\ast > 0$ such that $\mathbf{x}'=\mathbf{x}+t^\ast\mathbf{g}(C,\alpha)$. 
    The feasibility of $\mathbf{x}'$ implies that $\mathbf{x}_e<1$ for every $e \in C^+$ and $\mathbf{x}_e>0$ for every $e \in C^-$, for otherwise an entry of $\mathbf{x}'$ would be strictly bigger than one or negative. However, since $\mathbf{x}$ is the incidence vector of $M$, this implies that $C^+=C \setminus M$ and $C^-=C\cap M$. Thus $C$ is an alternating cycle for $M$, and $M':=M \Delta C$ is another perfect matching of $M$. By Lemma~\ref{lemma:adjacency} $\mathbf{1}_{M'}$ is adjacent to $\mathbf{x}=\mathbf{1}_M$ on the skeleton of $P_G$. We now show that $\mathbf{x}':=\mathbf{1}_{M'}$, which will conclude the proof of the claim. To see this, note that by definition of a circuit-move, $t^\ast \ge 0$ was chosen maximal such that $\mathbf{x}+t^\ast\mathbf{g}(C,\alpha) \in P_G$. Since moving in the direction of the circuit $\mathbf{g}(C,\alpha)$ preserves all equality-constraints $(1)$ in the hyperplane description of $P_G$, the maximality implies that moving $t$ above $t^\ast$ would violate one of the inequality-constraints $(2)$. Thus, there has to exist $e \in C$ such that $\mathbf{x}'_e =0$. On the other hand, we know that $$0=\mathbf{x}'_e=\underbrace{\mathbf{x}_e}_{\in \{0,1\}}\pm \underbrace{t^\ast\alpha}_{>0},$$ which implies that $t^\ast\alpha=1$. It follows that  $\mathbf{x}'=\mathbf{x}+t^\ast\mathbf{g}(C,\alpha)=\mathbf{1}_{M'}$, as desired. 
\end{enumerate}
\qed\end{proof}

\section*{Acknowledgements}

This work was initiated during a stay of the first author as a guest professor at ETH Zurich in May-July 2022. He wishes to thank Prof. Emo Welzl for his hospitality during this stay.

\bibliographystyle{plain}
\bibliography{pm.bib}

\end{document}